\documentclass[titlepage,twoside,12pt]{article}
\usepackage{amssymb}
\usepackage{amsfonts}
\begin{document}

{\bf \Large On the Safe Use of \\ \\ Inconsistent Mathematics} \\ \\

{\bf Elem\'{e}r E Rosinger} \\ \\
Department of Mathematics \\
and Applied Mathematics \\
University of Pretoria \\
Pretoria \\
0002 South Africa \\
eerosinger@hotmail.com \\ \\

{\bf Abstract} \\

A method is presented for using the consistent part of inconsistent axiomatic systems. \\ \\

{\bf 1. Half a Century of Practice of Inconsistent Mathematics} \\

Rather missed in the conscious awareness of many, we have been doing inconsistent mathematics for more than half a century
by now, and in fact, have quite heavily and essentially depended on it in our everyday life. Indeed, electronic digital
computers, when considered operating on integers, which is but a part of their operations, act according to the system of
axioms given by

\begin{itemize}

\item (PA) : the usual Peano Axioms for $\mathbb{N}$,

\end{itemize}

plus the {\it ad-hock axiom}, according to which

\begin{itemize}

\item (MI) : there exists $M \in \mathbb{N},~ M >> 1$, such that $ M + 1 = M$

\end{itemize}

Such an $M$, called "machine infinity", is usually larger than $10^{100}$, however, it is inevitably inherent in every
electronic digital computer, due to obvious unavoidable physical limitations. And clearly, the above mix of (PA) + (MI)
axioms is {\it inconsistent}. Yet we do not mind flying on planes designed and built with the use of such electronic
digital computers. \\ \\

{\bf 2. A Way to Handle General Inconsistencies} \\

The example of the mentioned quite elementary inconsistency characterizing each and every electronic digital computer can
suggest a way to deal with rather arbitrary inconsistent axiomatic systems. \\

Namely, the way inconsistencies in computation with integers are avoided in electronic digital computers is simply by {\it never reaching overflow}. In
other words, it is made sure that additions, subtractions and multiplications lead to results well within the range of integers between
$- M$ and $M$. It is of course more difficult to avoid overflow when operating on decimal numbers, since in such a case
one also has to guard against divisions with numbers that are too small. \\

However, such a rather simple manner of avoiding inconsistencies in a system which is essentially inconsistent can be
extrapolated far beyond the above case of (PA) + (MI) axioms. Here we shall indicate such a way, which may be seen as a simplification of an idea of
G\"{o}del, and pursued later by Ehrenfeucht \& Mycielski, see Manin [pp. 258-260]. \\

Let ${\cal L}$ be a {\it formal language} over an {\it alphabet} $A$. Let ${\cal W}$ be a set of {\it well formed formulas}
in ${\cal L}$. Lastly, let ${\cal D}$ be a {\it deduction rule} which for every subset of well formed formulas ${\cal S}
\subseteq {\cal W}$ gives the subset of well formed formulas ${\cal D} ( {\cal S} ) \subseteq {\cal W}$ which are the
consequences of ${\cal S}$. We shall later specify a certain finitary type condition on this mapping \\

(2.1)$~~~~~~ {\cal D} : {\cal P} ( {\cal W} ) \longrightarrow {\cal P} ( {\cal W} )$ \\

where as usual, ${\cal P} ( {\cal W} )$ denotes the set of all subsets of ${\cal W}$. \\

A {\it theory}\, ${\cal T} \subseteq {\cal W}$\, in such a setup $( A, {\cal L}, {\cal W}, {\cal D} )$ can be seen as the
result of the following construction. Given a set of {\it axioms} ${\cal A} \subseteq {\cal W}$, one obtains ${\cal T}$ by
the sequence of operations \\

(2.2)$~~~~~~ \begin{array}{l}
            {\cal T}_0 ( {\cal A} ) = {\cal A} \\
            {\cal T}_1 ( {\cal A} ) = {\cal D} ( {\cal T}_0 ( {\cal A} ) ) \\
            {\cal T}_2 ( {\cal A} ) = {\cal D} ( {\cal T}_1 ( {\cal A} ) ) \\
            \vdots \\
            {\cal T}_{n+1} ( {\cal A} ) = {\cal D} ( {\cal T}_n ( {\cal A} ) ) \\
            \vdots
         \end{array} $ \\

according to \\

(2.3)$~~~~~~ {\cal T} = {\cal T} ( {\cal A} ) = \bigcup_{n \in \mathbb{N}}~ {\cal T}_n ( {\cal A} ) $ \\

In the case of an inconsistent set ${\cal A}$ of axioms, one typically has \\

(2.4)$~~~~~~ {\cal T} ( {\cal A} )  = {\cal W} $ \\

and in the sequel, we shall consider (2.4) as the definition of the {\it inconsistency} of the set ${\cal A}$ of
axioms. \\

Let us now give a way to avoid inconsistencies in a theory ${\cal T} ( {\cal A} )$ which satisfies (2.4), thus itself
as a whole is inconsistent. \\
Here by avoiding inconsistency we mean identifying certain sentences $P \in {\cal T} ( {\cal A} )$ for which exists a subset of axioms ${\cal B}
\subseteq {\cal A}$, such that $P \in {\cal T} ( {\cal B} )$, while ${\cal B}$ is not inconsistent, that is, we have ${\cal T} ( {\cal B} ) \subsetneqq
{\cal W}$. \\

For every $P \in {\cal T} ( {\cal A} ) $, and in view of (2.2), (2.3), we denote by \\

(2.5)$~~~~~~ n_{{\cal A}, P} \in \mathbb{N} $ \\

the smallest $n \in \mathbb{N}$, such that $P \in {\cal T}_n ( {\cal A} )$. Obviously, $n_{{\cal A}, P}$ is well
defined. \\

Let us now assume that the mapping (2.1) is {\it finitely generated}, that is, it satisfies the condition \\

(2.6)$~~~~~~ {\cal D} ( {\cal S} ) = \bigcup {\cal D} ( {\cal F} ) $ \\

for all ${\cal S} \subseteq {\cal W}$, where the union is taken over all finite subsets ${\cal F} \subseteq {\cal S}$. \\

Then it is easy to see that the following properties result \\

(2.7)$~~~~~~ {\cal D} ( {\cal S} ) \subseteq {\cal D} ( {\cal S}\,' ) $ \\

for ${\cal S} \subseteq {\cal S}\,' \subseteq {\cal W}$. \\

Also, for every family ${\cal S}_i \subseteq {\cal W}$, with $i \in I$, we have \\

(2.8)$~~~~~~ \bigcup_{i \in I} {\cal D} ( {\cal S}_i ) \subseteq {\cal D} ( \bigcup_{i \in I} {\cal S}_i ) $ \\

Furthermore, if the family ${\cal S}_i \subseteq {\cal W}$, with $i \in I$, is right directed with respect to inclusion, then \\

(2.9)$~~~~~~ \bigcup_{i \in I} {\cal D} ( {\cal S}_i ) = {\cal D} ( \bigcup_{i \in I} {\cal S}_i ) $ \\

Given now $P \in {\cal T}_n ( {\cal A} )$, for a certain $n \in \mathbb{N},~ n \geq 1$, then (2.2) leads to finite
subsets \\

(2.10)$~~~~~~ {\cal S}_{n-1} \subseteq {\cal T}_{n-1} ( {\cal A} ),~
             {\cal S}_{n-2} \subseteq {\cal T}_{n-2} ( {\cal A} ),~ \ldots ,~
             {\cal S}_0 \subseteq {\cal T}_0 ( {\cal A} ) = {\cal A} $ \\

such that \\

(2.11)$~~~~~~ P \in {\cal D} ( {\cal S}_{n-1} ),~ {\cal S}_{n-1} \subseteq {\cal D} ( {\cal S}_{n-2} ),~ \ldots ,~
             {\cal S}_1 \subseteq {\cal D} ( {\cal S}_0 ) $ \\

Indeed, the sequence ${\cal S}_{n-1},~ {\cal S}_{n-2},~ \ldots ,~ {\cal S}_0$ in (2.10), (2.11) results from (2.2), and thus we only have to show that
its terms can be chosen as finite subsets. Here we note that (2.6), (2.11) give a finite subset ${\cal F}_{n-1} \subseteq {\cal S}_{n-1}$, such that $P
\in {\cal D} ( {\cal F}_{n-1} )$. Further, for every $Q_{n-1} \in {\cal F}_{n-1}$, there is a finite subset ${\cal F}_{Q_{n-1}} \subseteq {\cal
S}_{n-2}$, such that $Q_{n-1} \in {\cal D} ( {\cal F}_{Q_{n-1}} )$. Thus in view of (2.8), we have \\

$~~~~~~ {\cal F}_{n-1} \subseteq \bigcup_{\,Q_{n-1} \in\, {\cal F}_{n-1}}\, {\cal D} ( {\cal F}_{Q_{n-1}} ) \subseteq
                      {\cal D} ( \bigcup_{\,Q_{n-1} \in\, {\cal F}_{n-1}}\, {\cal F}_{Q_{n-1}} ) $ \\

and obviously, $\bigcup_{\,Q_{n-1} \in\, {\cal F}_{n-1}}\, {\cal F}_{Q_{n-1}}$ is a finite subset. And the argument follows in the same way, till it
reaches ${\cal S}_0$. \\

We denote by \\

(2.12)$~~~~~~ \Sigma_{{\cal A}, P} $ \\

the set of all such sequences ${\cal S}_{n-1},~ {\cal S}_{n-2},~ \ldots ,~ {\cal S}_0$ which satisfy (2.10), (2.11). \\

The interest in $\Sigma_{{\cal A}, P}$ is obvious. Indeed, if \\

(2.13)$~~~~~~ \exists~~ ( {\cal S}_{n-1},~ {\cal S}_{n-2},~ \ldots ,~ {\cal S}_0 ) \in \Sigma_{{\cal A}, P} ~:~
                                   {\cal T} ( {\cal S}_0 ) \subsetneqq {\cal W} $ \\

then in view of (2.10), (2.11), the sentence $P \in  {\cal T}_n ( {\cal A} )$ belongs to a {\it consistent} part of the theory ${\cal T} ( {\cal A} )$.
Therefore, we denote by \\

(2.14)$~~~~~~ {\cal C}_{\cal A} $ \\

the set of all such sentences $P$. \\

Finally, a rather simple way to identify sentences in ${\cal C}_{\cal A}$, a way suggested by the usual method of avoiding "machine infinity", is as
follows. \\

For every $P \in {\cal T}_n ( {\cal A} )$, with $n \in \mathbb{N},~ n \geq 1$, we define its {\t index} given by the two natural numbers \\

(2.15)$~~~~~~ \iota_{{\cal A}, P} = (  n_{{\cal A}, P},  m_{{\cal A}, P} ) $ \\

where $ n_{{\cal A}, P}$ was defined in (2.5), while for the second one we have \\

(2.16)$~~~~~~  m_{{\cal A}, P} = \min\, \sum\, length ( S ) $ \\

where \\

(2.17)$~~~~~~ length : {\cal W} \longrightarrow \mathbb{N} $ \\

and $length ( S )$, for $S \in {\cal W}$, denotes the number of letters in the alphabet $A$, with the possible repetitions counted, which constitute $S$,
while the minimum is taken over all sequences  $( {\cal S}_{n-1},~ {\cal S}_{n-2},~ \ldots ,~ {\cal S}_0 ) \in \Sigma_{{\cal A}, P}$, and the sum is
taken over all $S \in {\cal S}_{n-1} \cup {\cal S}_{n-2} \cup \ldots \cup {\cal S}_0$. \\

We define now \\

(2.18)$~~~~~~ \kappa_{\cal A} = \min\, \iota_{{\cal A}, P} $ \\

where the minimum is taken in the lexicographic order on $\mathbb{N} \times \mathbb{N}$, and over all $P \in {\cal T} ( {\cal A} ) \setminus
{\cal C}_{\cal A}$. \\

Then in view of the above, one obtains \\

{\bf Theorem 2.1.} \\

If for $P \in {\cal T} ( {\cal A} )$ one has in the lexicographic order on $\mathbb{N} \times \mathbb{N}$ \\

(2.19)$~~~~~~ \iota_{{\cal A}, P} < \kappa_{\cal A} $ \\

then \\

(2.20)$~~~~~~ P \in {\cal C}_{\cal A} $ \\

{\bf Remark 2.1.} \\

1) An upper estimate of $\kappa_{\cal A}$ can easily be obtained, since it is enough to make use of one single $P \in {\cal T} ( {\cal A} ) \setminus
{\cal C}_{\cal A}$. \\

2) One can obviously replace the axioms ${\cal A} \subseteq {\cal W}$ with an equivalent set of axioms ${\cal A} \subseteq {\cal W}$, which means that \\

(2.21)$~~~~~~ {\cal T} ( {\cal A} ) = {\cal T} ( {\cal B} ) $ \\

And as follows under appropriate conditions from Ehrenfeucht \& Mycielski, see also Manin [pp. 258-260], with such a replacement one may rather
arbitrarily shorten, or for that matter, lengthen the proofs of sentences $P \in {\cal T} ( {\cal A} ) = {\cal T} ( {\cal B} )$. Such an equivalent
replacement of axioms, therefore, can rather arbitrarily change $\iota_{{\cal A}, P}$, and hence also $\kappa_{\cal A}$. And obviously, any increase in
the latter makes the result in Theorem 2.1. above more strong, that is, it weakens the condition (2.19), and thus it increases the eventuality of
condition (2.20). \\

3) In the above, the mapping (2.1) was not required to satisfy the condition \\

(2.22)$~~~~~~ {\cal S} \subseteq {\cal D} ( {\cal S} ) $ \\

for ${\cal S} \subseteq {\cal W}$. \\

4) The more general use of inconsistent mathematics, [4], than merely in our electronic digital computers, as well as of mathematics developed upon
self-referential logic, [1], could become future trends. Both obviously open ways to realms not studied or used so far, and the latter brings with it
impressive possibilities in dealing with structures with complexities not encountered so far. \\ \\

{\bf 3. Comments} \\

A brief consideration of the millennia old relevant historical background may be appropriate. \\

Among fundamental ideas that have for long ages preoccupied humans, and apparently did so even prior to the existence of literate cultures and
civilizations, have been those of self-reference, infinity and change. These ideas appear time and again in a variety of forms of expression, and as
such, the wonder and puzzlement they have always produced and still keep producing have not been lost until our own days. \\

Self-reference, for instance, can be found illustrated since prehistoric times by a snake which bites its own tail. It is also clearly expressed in the
ancient Greek Paradox of the Liar which, in its modern version and within set theory, takes the form of Russell's Paradox. \\
Needless to say, with such a record, the idea took quite strongly hold ever since ancient times that one should better avoid self-reference, since it
does so easily lead to paradoxes. And the effect of that persistent and pervasive view has been the association of the manifestly negative term of
"vicious circle" with any form of self-reference. \\
On the other hand, and apparently prior to the emergence of the Paradox of the Liar, the Old Testament, in Exodus 3:14, gives the name of God by what is
for all purposes an ultimate and all encompassing self-reference, namely, "I Am that I Am". \\

Infinity and change is, among others, the substance of Zeno's surprising paradoxes which until our own days have not been fully clarified. Yet, unlike
with self-reference, neither infinity, nor change got such a negative aura as to impel us into their avoidance. And in fact, since Newton's Calculus,
both change and infinity have become fundamental in Mathematics, with their extraordinary consequences for technology. \\

What has, however, appeared to had been perfectly clarified since time immemorial is logical contradiction. Indeed, there does not appear to be any
record of an instance when, due to some reason, the decision was taken to avoid contradiction at all costs, and do so without any exception. After all,
paradoxes themselves are paradoxes only if we outlaw logical contradiction, and do so in absolutely every case when it may appear. Thus in particular,
the Paradoxes of the Liar or of Zeno would instantly lose much of their significance if we were to be more lenient with logical contradictions. \\

It is against such a background, therefore, that "inconsistent mathematics" and the "vicious circles" of "non-well-founded" set theory come to the
fore. \\

However, several aspects in this regard are worth noting. \\

First perhaps, is the fact that, at least since the mid 1940s, when we started more and more massively using and relying upon electronic digital
computers, we have in fact been permanently practicing "inconsistent mathematics", including in may extremely important situations, and were doing so
without giving it any special thought. In this regard it may be noted that Mortensen's 1995 book itself gives hardly any consideration to the obvious and
elementary logical contradiction between the Peano Axioms and the inevitable existence of a so called "machine infinity" in each and every electronic
digital computer. \\

The good news in this regard, nevertheless, is that at last we do no longer keep being lost in that immemorial decision to avoid logical contradiction
by all means and without any exception. \\

As for the essentially self-referential nature of the so called "non-well-founded" set theory in the 1996 book of Barwise and Moss, its original
motivation is to a good extent also related to electronic digital computers, namely, to the construction of the various formal languages which facilitate
their use. \\
In this way, our long not noted or disregarded essential involvement in "inconsistent mathematics" came about as an imposition by the hardware of
electronic digital computers. On the other hand, our relatively recent interest in the self-referential "vicious circles" of "non-well-founded" set
theory was inspired by important software aspects related to electronic digital computers. \\

For the time being, however, each of these developments seem to proceed strictly on its own way. Indeed, "inconsistent mathematics" is not getting
involved in self-reference. As for the self-referential "vicious circles" of "non-well-founded" set theory, it is very carefully kept away from any sort
of logical contradictions, that accomplishment being in fact one of its basic aims, as well as claims to success. \\

Needless to say, both these developments can turn into extraordinary openings in Mathematics and in Mathematical Logic. And yet, a far greater opening is
there, in bringing them together at some future time. \\
Until then, by avoiding any sort of logical contradictions, the "vicious circles" of "non-well-founded" set theory can develop sufficiently and gain a
well deserved respectability. After that, the inclusion of logical contradictions in its further development may contribute significantly to its further
development. \\

And such a development may simply mean no more than setting at last aside two ancient and deeply rooted habits of thinking, namely, horror of logical
contradictions, and extreme suspicion of self-reference. \\

\end{document}